\documentclass[12pt]{article}

\usepackage{pstricks,pst-node}
\usepackage{psfig}
\usepackage{amsmath,amsthm}
\usepackage{amssymb}
\usepackage{color}
\usepackage{xspace}
\usepackage[colorlinks=true,
linkcolor=green,
filecolor=brown,
citecolor=green]{hyperref}

\def\v{\vert}
\def\eps{\epsilon}

\def\s{\ensuremath{\mathcal S}\xspace}
\def\l{\ensuremath{\mathcal L}\xspace}
\def\k{\ensuremath{\mathcal K}\xspace}
\def\p{\ensuremath{\mathcal P}\xspace}
\def\r{\ensuremath{\mathcal R}\xspace}

\begin{document}
\newtheorem{theorem}{Theorem}
\newtheorem{defn}[theorem]{Definition}
\newtheorem{lemma}[theorem]{Lemma}
\newtheorem{prop}[theorem]{Proposition}
\newtheorem{cor}[theorem]{Corollary}
\begin{center}
{\Large
Some bijections for lattice paths                        \\ 
}

\vspace{2mm}
David Callan  \\
\vspace{5mm}
Dec 9, 2021
\end{center}

\begin{abstract}
We present three bijections, the first between little Schr\"{o}der paths and a class of growth-constrained integer sequences, the second between lattice paths consisting of steps with nonnegative slope and another class of growth-constrained sequences, the third between a class of lattice paths with steps $(1,1)$ and $(1,-j), \ j\ge 1,$ and a class with steps $(k,1),\ k\ge 1,$ and $(1,-1)$.
\end{abstract}

\section{Introduction} \label{intro}

\vspace*{-4mm}

The purpose of this paper is to present three bijections involving combinatorial objects counted by sequences 
\htmladdnormallink{A001003}{http://oeis.org/A001003}, 
\htmladdnormallink{A049600}{http://oeis.org/A049600}, and
\htmladdnormallink{A090345}{http://oeis.org/A090345} respectively
in The On-Line Encyclopedia of Integer Sequences (OEIS) \cite{oeis}. Lattice paths are a common feature and the bijections are presented in the next 3 sections.

\section{A manifestation of the little Schr\"{o}der numbers} \label{schr}

\vspace*{-4mm}

The little Schr\"{o}der numbers $(s_n)_{n\ge 1} =(1,3,11,45,197,\dots)$ form sequence 
\htmladdnormallink{A001003}{http://oeis.org/A001003} in OEIS. They have many combinatorial manifestations (see e.g., \cite[Sec. 6.2.8 and Ex. 6.39]{ec2}) including little  Schr\"{o}der paths, defined below. In this section we give a bijective proof for another one: $s_n=\v \s_n \v$ where $\s_n$ is the set of sequences of nonnegative integers $(u_1,u_2,\dots,u_n)$ such that (i) $u_1=1$, (ii) $u_i \le i$ for all $i$, and (iii) the nonzero $u_i$ are weakly increasing. For instance, $\s_2=\{10,11,12\}$.

A \emph{Schr\"{o}der path} of semilength $n$ is a lattice path of steps east $E=(1,0)$, north $N=(0,1)$ and northeast (diagonal) $D=(1,1)$ from the origin to $(n,n)$ that never rises above the diagonal line $y=x$. Such a path has equal numbers of $E$ and $N$ steps and its semilength is the number of $E$ steps plus the number of $D$ steps. A \emph{little} Schr\"{o}der path is one that has no $D$ steps on the diagonal $y=x$. Let $\l\s$ denote them.  For instance, $\l\s_2=\{EENN,EDN,ENEN\}$.

Here is a bijection from $\l\s_n$ to $\s_n$ for $n\ge 1$.
Given a little Schr\"{o}der path, say $EDENEDNDDNEN \in \l\s_8$ as in Figure 1 in blue, record the heights above the horizontal line $y=-1$ of the terminal points of each $E$ and $D$ step, marking the height of each diagonal step with a superscript $^d$ as in $h^d$.

\begin{center}
\begin{pspicture}(-5,-2.8)(4,8.5)

\psline[linecolor=blue](-4,0)(-3,0)(-2,1)(-1,1)(-1,2)(0,2)(1,3)(1,4)(2,5)(3,6)(3,7)(4,7)(4,8)
 
\psdots(-4,0)(-3,0)(-2,1)(-1,1)(-1,2)(0,2)(1,3)(1,4)(2,5)(3,6)(3,7)(4,7)(4,8)

\psset{linecolor=lightgray}
\psline(-4,-1)(4,-1)
\psline(-3,0)(4,0)
\psline(-1,1)(4,1)
\psline(0,2)(4,2)
\psline(1,3)(4,3)
\psline(1,4)(4,4)
\psline(2,5)(4,5)
\psline(3,6)(4,6)
\psline(4,-1)(4,7)
\psline(3,-1)(3,6)
\psline(2,-1)(2,5)
\psline(1,-1)(1,3)
\psline(0,-1)(0,2)
\psline(-1,-1)(-1,1)
\psline(-2,-1)(-2,1)
\psline(-3,-1)(-3,0)
\psline(-4,-1)(-4,0)
\rput(5.35,0){\textrm{{\small $x$-axis}}}
\rput(4.4,-0.05){\textrm{{\small $\leftarrow$}}}
\rput(4.4,-1){\textrm{{\small $\leftarrow$}}}
\rput(5.4,-1){\textrm{{\small $y=-1$}}}
\psdots(-3,1)(-2,2)(-1,3)(0,3)(0,4)(1,5)(2,6)
\psdots(-4,-1)(-3,-1)(-2,-1)(-1,-1)(0,-1)(1,-1)(2,-1)(3,-1)(4,-1)  

\rput(-4.5,.1){\textrm{{\footnotesize (0,0)}}}
\rput(4.5,8.2){\textrm{{\footnotesize (8,8)}}}

\rput(-5.7,-.6){\textrm{{\small height sequence\,:}}}
\rput(-3.5,-.6){\textrm{{\small 1}}}
\rput(-2.5,-.55){\textrm{{\small $2^d$}}}
\rput(-1.5,-.6){\textrm{{\small 2}}}
\rput(-0.5,-.6){\textrm{{\small $3$}}}
\rput(0.5,-.55){\textrm{{\small $4^d$}}}
\rput(1.5,-.55){\textrm{{\small $6^d$}}}
\rput(2.5,-.55){\textrm{{\small $7^d$}}}
\rput(3.5,-.6){\textrm{{\small 8}}}

\rput(0,-1.8){\textrm{{\small A little Schr\"{o}der path and its height sequence}}}
\rput(0,-2.7){\textrm{{\small Figure 1}}}  

\psset{linecolor=gray}
\psline(-3,0)(4,0)

\end{pspicture}
\end{center} 

The example produces the height sequence $1\,2^d\,2\,3\,4^d\,6^d\,7^d\,8$. 
Clearly, this height sequence determines the path and
a valid height sequence 
begins with an unmarked 1, is weakly increasing and strictly so if the second height is marked, and the $i$-th height is $\le i$ for all $i$.

Now, say a marked height $h^d$ is \emph{lonely} if its predecessor height, 
$g$ or $g^d$, satisfies $h \ge g+2$ (this occurs when a $D$ in the path is immediately preceded by an $N$). The main step in the bijection is iterative: as long as a lonely $h^d$ is present, interchange the first such $h^d$ (although any such $h^d$ would do) with its predecessor and decrement $h$ to $h-1$. The example yields
\[
\begin{array}{ccccccccr}
1 & 2^d &  2  &  3  & 4^d & 6^d & 7^d &  8   & \rightarrow \\
  &     &     &     & 5^d & 4^d &     &      & \rightarrow \\
  &     &     & 4^d & 3   &     &     &      & \rightarrow \\
  &     & 3^d &  2  &     &     &     &      & \rightarrow \\
1 & 2^d & 3^d &  2  & 3   & 6^d & 4^d &  8   & \rightarrow \\ 
  &     &     &     & 5^d &  3  &     &      & \rightarrow \\ 
1 & 2^d & 3^d & 4^d & 2   &  3  & 4^d &  8\,,&     
\end{array}
\]
where blank spaces indicate no change from the previous line.
When done, each marked height $h^d$ will be preceded by $(h-1)^d$ or $h-1$. The last step is to replace each marked height with 0, yielding
$1\,0\,0\,0\,2\,3\,0\,8$ in our example.

The iterative steps in the bijection trade ``lonely''  marked heights for ``offending''  marked heights, where $h^d$ is offending if it has a right neighbor $j$ or $j^d$ with $h>j$. It should now be clear how to reverse the map by eliminating offending marked heights.

In case there are no diagonal steps in a little Schr\"{o}der path, we have a (rotated) 
Dyck path, and the above bijection produces a well known manifestation of the Catalan numbers: sequences of positive integers $(u_1,u_2,\dots,u_n)$ such that
$u_i \le i$ for all $i$ and the $u_i$ are weakly increasing.

\section{Kimberling paths} \label{kimb}

\vspace*{-4mm}

A \emph{Kimberling path}, after \htmladdnormallink{A049600}{http://oeis.org/A049600}, is a lattice path from the origin $O=(0,0)$ to a point $(i,j)$ that consists of steps with finite nonnegative slope. Thus vertical steps (with infinite slope) are not allowed. Let $\k(i,j)$ denote the Kimberling paths that terminate at $(i,j)$ (such paths exist only for $i\ge 1$ and $j\ge 0$). The 3 paths in $\k(2,1)$ are shown in Figure 2. 
\begin{center}
\begin{pspicture}(-5,-2)(5,1.5)
\psline[linecolor=blue](-6,0)(-5,0)(-4,1)
\psline[linecolor=blue](-1,0)(0,1)(1,1)
\psline[linecolor=blue](4,0)(6,1)
 
\psdots(-5,0)(-6,0)(-4,1)(-1,0)(0,1)(1,1)(4,0)(6,1)
\psset{linecolor=lightgray}
\psdots(-4,0)(0,0)(1,0)(5,0)(6,0)
\psdots(-6,1)(-5,1)(-1,1)(4,1)(5,1)

\rput(-6,-.3){\textrm{{\footnotesize $O$}}}
\rput(-3.5,1.1){\textrm{{\footnotesize (2,1)}}}
\rput(-1,-.3){\textrm{{\footnotesize $O$}}}
\rput(1.5,1.1){\textrm{{\footnotesize (2,1)}}}
\rput(4,-.3){\textrm{{\footnotesize $O$}}}
\rput(6.5,1.1){\textrm{{\footnotesize (2,1)}}}

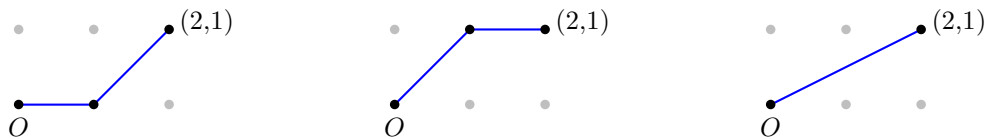
\rput(0,-1.2){\textrm{{\small The paths in $\k(2,1)$: the middle path has vertex set $\{(0,0),(1,1),(2,1)\}$}}} 

\rput(0,-2){\textrm{{\small Figure 2}}}  

\end{pspicture}
\end{center}  
Choosing the $x$-coordinates of the interior vertices of a Kimberling path and then the $y$-coordinates, the number of paths in $\k(i,j)$ with $k$ interior vertices is easily seen to be $\binom{i-1}{k}\binom{j+k}{k}$.

For $i,j \ge 0$, let $\l(i,j)$ denote the set of sequences of positive integers $(u_1,u_2,\dots,u_i)$ such that (i) $u_k \le j+2$ for all $k$, and (ii) $u_k\ge 
\max\{u_1,\dots,u_{k-1}\}-1$ for $k\ge 2$. By convention, $\l(0,j)$ consists of just one sequence, the empty sequence $\eps$ of length 0, and for example, $\l(1,j)=\{1,2,\dots,j+2\}$ and $\l(2,1)=\{11, 12, 13, 21, 22, 23, 32, 33\}$, where each sequence is presented as a contiguous string of digits.

The set  $\l(i,j)$ is equinumerous with $\k(i+1,j)$ and we now
define a bijection $\phi_j \,:\l(i,j) \to \k(i+1,j)$.
The subscript $j$ on $\phi$ is needed because, while the value of $i$ can be inferred from an integer sequence in $\l$ (being its length), the value of $j$ cannot. The bijection is defined recursively by induction on both $i$ and $j$. For brevity, we write $\phi_j(u_1,\dots,u_i)$ for $\phi_j(u)$ with $u=(u_1,\dots,u_i)\in \l(i,j)$.

For the base cases, with $i=0$, $\phi_j(\eps)$ is the path of one step from (0,0) to $(1,j)$, and with $i=1$, 
\[
\phi_j( 1 )=\{(0,0),(2,j)\} \textrm{ while } \phi_j(k)=\{(0,0),(1,k-2),(2,j)\} \textrm{ for } 2\le k \le j+2\, .
\]
Now suppose $i\ge 2$ and $u_1=k$. Then $\phi_j(k,u_2,\dots,u_i)$ is defined as follows.
If $k\le 2$, add (1,0) to each vertex of $\phi_j(u_2,\dots,u_i)$, that is, shift each vertex one unit east and then 
\begin{itemize}
\vspace*{-3mm}
\item if $k=1$, replace vertex (1,0) by (0,0) so the image path has no vertex on the vertical line $x=1$,
\item if $k=2$, prepend (0,0) so that the image path begins $\{(0,0),(1,0),\dots\}$. 
\end{itemize}
\vspace*{-3mm}
On the other hand, if $k\ge 3$, 
add $(1,k-2)$ to each vertex of $\phi_{j-(k-2)}(u_2-(k-2),\dots,u_i-(k-2))$ and then prepend $(0,0)$ so that $(1,k-2)$ is a vertex on the path for $3 \le k \le j+2$. For example, 
$\phi_3(1,4,3,5)=\{ (0,0),(2,2),(4,3),(5,3) \}$ because
$\phi_3(4,3,5)=\{ (0,0),(1,2),(3,3),(4,3)\}$ which follows from 
$\phi_1(1,3)=\{(0,0),(2,1),(3,1)\}$.
 
It is not hard to recursively define an inverse for $\phi$, looking first at where a Kimberling path intersects the line $x=1$. 
 
\section{Deutsch paths and Ram\'{i}rez paths} \label{deut}

\vspace*{-4mm}

A \emph{Deutsch path} \cite{DeutschPaths} is a lattice path of upsteps (1,1) and downsteps $(1,-j),\ j\ge 1$, that starts at the origin and never goes below the $x$-axis, and it is closed if it ends on the $x$-axis. Let $\p_n$ denote the set of $n$-step closed Deutsch paths with short valley downsteps, that is, each downstep followed by an upstep thereby creating a valley is of the from $(1,-1)$, the shortest possible. 
For example, $\p_4$ consists of $UUU3,\, UU11, \,U1U1$, where an 
upstep is denoted by $U$ and a downstep $(1,-j)$ is denoted by $j$ so that a number immediately preceding a $U$ must be 1, and $\p_5$  consists of $UUUU4,\, UUU12,\, UUU21,\, UU1U2,\, U1UU2$.
 
A \emph{Ram\'{i}rez path}, after  
\htmladdnormallink{A090345}{http://oeis.org/A090345}, 
is a lattice path of upsteps $(k,1)$ with $k\ge 1$  
and downsteps $(1,-1)$ that starts at the origin, never goes below the $x$-axis, 
and ends on the $x$-axis. 
Its size is the $x$-coordinate of its terminal point. Let $\r_n$ denote 
the set of Ram\'{i}rez paths of size $n$. For example, $\r_4$ consists of 
$1D1D, \,11DD,\, 3D$ where a downstep is denoted by $D$ and an upstep $(k,1)$ 
is denoted by $k$, and $\r_5$ consists of $1D2D,\, 12DD,\, 2D1D,\, 21DD,\,4D$.

We will show that $\p_n$ and $\r_n$ are equinumerous. 
Note that $\p_0$ and $\r_0$ both consist of the empty path, 
$\p_1$ and $\r_1$ are both empty, and $\p_2$ and $\r_2$ both consist 
of the single path $UD$.
We now give a simple graphical description of a bijection from $\p_n$ to $\r_n$ 
for $n\ge 3$. So suppose given a path $P\in 
\p_n$ as in Figure 3a) below. The image path $Q \in \r_n$ is obtained as 
follows. First replace each ``long'' downstep $(1,-j), \ j\ge 2$ with $j$ short 
downsteps $(1/j,-1)$ and color the last $j-1$ of them blue along with their matching 
upsteps as in Figure 3b). 
 
\begin{center}
\begin{pspicture}(-1,-4)(28,3)
\psset{unit=.5cm}
\psline[linewidth=.05cm](-1,-1)(0,0)(1,1)(2,2)(3,3)(4,2)(5,3)(6,4)(7,5)(8,4)(9,5)(10,6)(11,3)(12,2)(13,-1)

\psline[linecolor=gray](-1,-1)(13,-1)
\psline[linecolor=gray](16,-1)(30,-1)

\psset{dotsize=6pt 0}
\psdots(-1,-1)(0,0)(1,1)(2,2)(3,3)(4,2)(5,3)(6,4)(7,5)(8,4)(9,5)(10,6)(11,3)(12,2)(13,-1)
\psset{dotsize=5pt 0}
\psset{linecolor=lightgray}
\psdots(-1,-1)(0,-1)(1,-1)(2,-1)(3,-1)(4,-1)(5,-1)(6,-1)(7,-1)(8,-1)(9,-1)(10,-1)(11,-1)(12,-1)(13,-1) 

\rput(-1,-1.5){\textrm{{\scriptsize 0}}}
\rput(0,-1.5){\textrm{{\scriptsize 1}}}
\rput(1,-1.5){\textrm{{\scriptsize 2}}}
\rput(2,-1.5){\textrm{{\scriptsize 3}}}
\rput(3,-1.5){\textrm{{\scriptsize 4}}}
\rput(4,-1.5){\textrm{{\scriptsize 5}}}
\rput(5,-1.5){\textrm{{\scriptsize 6}}}
\rput(6,-1.5){\textrm{{\scriptsize 7}}}
\rput(7,-1.5){\textrm{{\scriptsize 8}}}
\rput(8,-1.5){\textrm{{\scriptsize 9}}}
\rput(9,-1.5){\textrm{{\scriptsize 10}}}
\rput(10,-1.5){\textrm{{\scriptsize 11}}}
\rput(11,-1.5){\textrm{{\scriptsize 12}}}
\rput(12,-1.5){\textrm{{\scriptsize 13}}}
\rput(13,-1.5){\textrm{{\scriptsize 14}}}
\rput(6,-2.8){\textrm{{\normalsize $P = U\,U\,U\,U\,1\,U\,U\,U\,1\,U\,U\,3\,1\,3 \in \p_{14}$}}}

\rput(6,-4){\textrm{{\normalsize A closed Deutsch path with }}} 
\rput(6,-5){\textrm{{\normalsize short valley downsteps}}} 
\rput(6,-6.5){\textrm{{\normalsize Figure 3a)}}} 
\rput(23,-6.5){\textrm{{\normalsize Figure 3b)}}} 

\psset{linecolor=black}

\rput(14.7,2){\textrm{\Large{$\longrightarrow$}}}  
\psline[linewidth=.05cm](16,-1)(17,0)(18,1)(19,2)(20,3)(21,2)(22,3)(23,4)(24,5)(25,4)(26,5)(27,6)(27.33,5)(27.66,4)(28,3)(29,2)(29.33,1)(29.66,0)(30,-1)

\psline[linecolor=blue,arrows=->,arrowsize=3pt 4](27.1,4.5)(26,4.5)
\psline[linecolor=blue,arrows=->,arrowsize=3pt 4](27.4,3.5)(23,3.5)
\psline[linecolor=blue,arrows=->,arrowsize=3pt 4](29.1,0.5)(18,0.5)
\psline[linecolor=blue,arrows=->,arrowsize=3pt 4](29.5,-0.5)(17,-0.5)

\psset{dotsize=6pt 0}
\psdots(16,-1)(17,0)(18,1)(19,2)(20,3)(21,2)(22,3)(23,4)(24,5)(25,4)(26,5)(27,6)(28,3)(29,2)(30,-1)
\psset{linewidth=.07cm,linecolor=blue}
\psline(27.33,5)(27.66,4)(28,3)
\psline(29.33,1)(29.66,0)(30,-1)
\psline(16,-1)(17,0)(18,1)
\psline(22,3)(23,4)
\psline(25,4)(26,5)
\psset{linewidth=.05cm}
\psset{dotsize=4pt 0}
\psdots(27.33,5)(27.66,4)(29.33,1)(29.66,0)(0,0)

\psset{dotsize=5pt 0}
\psset{linecolor=lightgray}
\psdots(16,-1)(17,-1)(18,-1)(19,-1)(20,-1)(21,-1)(22,-1)(23,-1)(24,-1)(25,-1)(26,-1)(27,-1)(28,-1)(29,-1)(30,-1) 
\rput(16,-1.5){\textrm{{\scriptsize 0}}}
\rput(17,-1.5){\textrm{{\scriptsize 1}}}
\rput(18,-1.5){\textrm{{\scriptsize 2}}}
\rput(19,-1.5){\textrm{{\scriptsize 3}}}
\rput(20,-1.5){\textrm{{\scriptsize 4}}}
\rput(21,-1.5){\textrm{{\scriptsize 5}}}
\rput(22,-1.5){\textrm{{\scriptsize 6}}}
\rput(23,-1.5){\textrm{{\scriptsize 7}}}
\rput(24,-1.5){\textrm{{\scriptsize 8}}}
\rput(25,-1.5){\textrm{{\scriptsize 9}}}
\rput(26,-1.5){\textrm{{\scriptsize 10}}}
\rput(27,-1.5){\textrm{{\scriptsize 11}}}
\rput(28,-1.5){\textrm{{\scriptsize 12}}}
\rput(29,-1.5){\textrm{{\scriptsize 13}}}
\rput(30,-1.5){\textrm{{\scriptsize 14}}}

\end{pspicture}

\begin{pspicture}(-.5,-3)(30,3.5)
\psset{unit=.5cm}
\psline[linewidth=.05cm](-1,-1)(0,0)(1,1)(2,2)(3,3)(4,2)(5,3)(6,4)(7,5)(8,4)(9,5)(10,6)(10.33,5)(10.66,4)(11,3)(12,2)(12.33,1)(12.66,0)(13,-1)
\psline[linecolor=gray](-1,-1)(13,-1)
\psline[linecolor=gray](16,-1)(30,-1)

\psline[linecolor=blue,arrows=->,arrowsize=3pt 4](10.1,4.5)(9,4.5)
\psline[linecolor=blue,arrows=->,arrowsize=3pt 4](10.4,3.5)(6,3.5)
\psline[linecolor=blue,arrows=->,arrowsize=3pt 4](12.1,0.5)(1,0.5)
\psline[linecolor=blue,arrows=->,arrowsize=3pt 4](12.5,-0.5)(0,-0.5)

\psset{dotsize=6pt 0}
\psdots(-1,-1)(0,0)(1,1)(2,2)(3,3)(4,2)(5,3)(6,4)(7,5)(8,4)(9,5)(10,6)(10.33,5)(10.66,4)(11,3)(12,2)(12.33,1)(12.66,0)(13,-1)
\psset{linewidth=.07cm,linecolor=blue}
\psline(10.33,5)(10.66,4)(11,3)
\psline(12.33,1)(12.66,0)(13,-1)
\psline(-1,-1)(0,0)(1,1)(2,2)
\psline(5,3)(6,4)(7,5)
\psline(8,4)(9,5)(10,6)
\psset{dotsize=4pt 0}
\psdots(10.33,5)(10.66,4)(12.33,1)(12.66,0)(0,0)(1,1)(6,4)(9,5)

\psset{linewidth=.05cm}
\psset{dotsize=5pt 0}
\psset{linecolor=lightgray}
\psdots(-1,-1)(0,-1)(1,-1)(2,-1)(3,-1)(4,-1)(5,-1)(6,-1)(7,-1)(8,-1)(9,-1)(10,-1)(11,-1)(12,-1)(13,-1)

\psset{linecolor=black}
\rput(15,1){\textrm{\Large{$\longrightarrow$}}}  
\psline[linewidth=.05cm](16,-1)(19,0)(20,1)(21,0)(22,1)(24,2)(25,1)(27,2)(28,1)(29,0)(30,-1)

\psset{dotsize=6pt 0}
\psdots(16,-1)(19,0)(20,1)(21,0)(22,1)(24,2)(25,1)(27,2)(28,1)(29,0)(30,-1)

\psset{dotsize=5pt 0}
\psset{linecolor=lightgray}
\psdots(16,-1)(17,-1)(18,-1)(19,-1)(20,-1)(21,-1)(22,-1)(23,-1)(24,-1)(25,-1)(26,-1)(27,-1)(28,-1)(29,-1)(30,-1)

\rput(-1,-1.5){\textrm{{\scriptsize 0}}}
\rput(0,-1.5){\textrm{{\scriptsize 1}}}
\rput(1,-1.5){\textrm{{\scriptsize 2}}}
\rput(2,-1.5){\textrm{{\scriptsize 3}}}
\rput(3,-1.5){\textrm{{\scriptsize 4}}}
\rput(4,-1.5){\textrm{{\scriptsize 5}}}
\rput(5,-1.5){\textrm{{\scriptsize 6}}}
\rput(6,-1.5){\textrm{{\scriptsize 7}}}
\rput(7,-1.5){\textrm{{\scriptsize 8}}}
\rput(8,-1.5){\textrm{{\scriptsize 9}}}
\rput(9,-1.5){\textrm{{\scriptsize 10}}}
\rput(10,-1.5){\textrm{{\scriptsize 11}}}
\rput(11,-1.5){\textrm{{\scriptsize 12}}}
\rput(12,-1.5){\textrm{{\scriptsize 13}}}
\rput(13,-1.5){\textrm{{\scriptsize 14}}}
\rput(16,-1.5){\textrm{{\scriptsize 0}}}
\rput(17,-1.5){\textrm{{\scriptsize 1}}}
\rput(18,-1.5){\textrm{{\scriptsize 2}}}
\rput(19,-1.5){\textrm{{\scriptsize 3}}}
\rput(20,-1.5){\textrm{{\scriptsize 4}}}
\rput(21,-1.5){\textrm{{\scriptsize 5}}}
\rput(22,-1.5){\textrm{{\scriptsize 6}}}
\rput(23,-1.5){\textrm{{\scriptsize 7}}}
\rput(24,-1.5){\textrm{{\scriptsize 8}}}
\rput(25,-1.5){\textrm{{\scriptsize 9}}}
\rput(26,-1.5){\textrm{{\scriptsize 10}}}
\rput(27,-1.5){\textrm{{\scriptsize 11}}}
\rput(28,-1.5){\textrm{{\scriptsize 12}}}
\rput(29,-1.5){\textrm{{\scriptsize 13}}}
\rput(30,-1.5){\textrm{{\scriptsize 14}}}

\rput(22.5,-2.8){\textrm{{\normalsize $Q=3\,1\,D\,1\,2\,D\,2\,D\,D\,D  \in \r_{14}$}}}

\rput(22.5,-4){\textrm{{\normalsize A Ram\'{i}rez path}}}  
\rput(6,-5.5){\textrm{{\normalsize Figure 4a)}}} 
\rput(22.5,-5.5){\textrm{{\normalsize Figure 4b)}}}

\end{pspicture}

\end{center} 
Each (maximal) run of blue upsteps in Figure 3b) is necessarily followed by an 
(uncolored) upstep; color this upstep blue also as in Figure 4a). Lastly, replace each 
run of $k$ blue upsteps by a step $(k,1)$, delete all the blue downsteps, and 
change $(1/j,-1)$ downsteps to regular $(1,-1)$ downsteps, as in Figure 4b).

We leave the reader to devise the inverse map.

\noindent Department of Statistics, University of Wisconsin-Madison \\
{\bf callan@stat.wisc.edu}

\end{document}